\documentclass[11pt]{amsart}
\usepackage{amsmath}
\usepackage{amsthm}

\newcommand{\Q}{{\mathbb Q}}

\newcommand{\Z}{{\mathbb Z}}
\newcommand{\C}{{\mathbb C}}

\renewcommand{\O}{{\mathcal O}}
\renewcommand{\H}{{\mathbb H}}

\newtheorem{theorem}{Theorem}
\newtheorem{lemma}{Lemma}
\newtheorem{corollary}{Corollary}
\newtheorem{proposition}[theorem]{Proposition}
\numberwithin{theorem}{subsection}
\numberwithin{corollary}{subsection}
\numberwithin{lemma}{subsection}
\theoremstyle{definition}
\newtheorem{conj}{Conjecture}
\numberwithin{conj}{subsection}
\newtheorem{example}{Example}
\numberwithin{example}{subsection}
\newtheorem{definition}{Definition}
\numberwithin{definition}{subsection}
\newtheorem{question}{Question}
\numberwithin{question}{subsection}
\numberwithin{equation}{subsection}
\theoremstyle{remark}
\newtheorem{remark}{Remark}
\numberwithin{remark}{subsection}
\begin{document}

\newcommand{\ra}{\rightarrow}

\title{Frobenius splitting and ordinarity}
\author{Kirti Joshi and C.~S.~Rajan}

\address{Department of
Mathematics, University of Arizona, 617 N Santa Rita, P O Box 210089,
Tucson, AZ 85721, USA, School of Mathematics, Tata Institute of Fundamental 
Research, Homi Bhabha Road, Bombay - 400 005, INDIA.}
\email{kirti@math.arizona.edu,rajan@math.tifr.res.in}


\let\tensor=\otimes 
\newcommand{\Hom}{{{\rm Hom}}}
\newcommand{\Spec}{{{\rm Spec}}}
\let\tensor=\otimes
\let\isom=\equiv 
\renewcommand{\P}{{\mathbb P}}
\let\liminv=\varprojlim

\begin{abstract}
   We examine the relationship between the notion of Frobenius
splitting and ordinarity for varieties. We show the following: a) The
de Rham-Witt cohomology groups $H^i(X, W({\mathcal O}_X))$ of a smooth
projective Frobenius split variety are finitely generated over $W(k)$.
b) we provide counterexamples to a question of V.~B.~Mehta that
Frobenius split varieties are ordinary or even Hodge-Witt. c) a Kummer
$K3$ surface associated to an Abelian surface is $F$-split (ordinary)
 if and only if the associated Abelian surface is $F$-split (ordinary). d) for a $K3$-surface
defined over a number field, there is a set of primes of density one
in some finite extension of the base field, over which the surface
acquires ordinary reduction.
\end{abstract}

\maketitle

\section{Introduction}
   Let $k$ be a perfect field of characteristic $p>0$. An abelian
variety $A$ over $k$ is said to be ordinary if the $p$-rank of $A$ is
the maximum possible, namely equal to the dimension of $A$. The notion
of ordinarity was extended by Mazur \cite{mazur72}, to a smooth
projective variety $X$ over $k$, using notions from crystalline
cohomology. A more general definition was given by Bloch-Kato
\cite{bloch86} and Illusie-Raynaud \cite{illusie83b} using coherent
cohomology.  With this definition it is easier to see that ordinarity
is an open condition. Ordinary varieties tend to have special
properties, for example the existence of canonical Serre-Tate liftings
for ordinary abelian varieties to characteristic $0$, and the
comparison theorems between crystalline cohomology and p-adic \'etale
cohomology can be more easily established for such varieties.  In
brief, ordinary varieties play a key role in the study of varieties in
characteristic $p>0$.

      One of the motivating questions in this paper is to study the
relationship between the concept of Frobenius split varieties
introduced by Mehta-Ramanathan \cite{mehta85b}, and ordinary varieties.
Unlike ordinariness, the definition of Frobenius splitting can be
extended to singular varieties, and this has proved to be of use in
studying the cohomology and singularities of Schubert varieties. It
was shown by Mehta and Srinivas \cite{mehta87}, that smooth,
projective varieties with trivial cotangent bundle, in particular
abelian varieties, are ordinary if and only if they are Frobenius
split. It can be seen from the theory of Cartier operator, that
smooth, projective $F$-split surfaces are ordinary (see
Theorem~\ref{surface-fsord}). Moreover any smooth, projective,
ordinary variety with trivial canonical bundle is $F$-split (see
Theorem~\ref{ktrivial-ordfs}).

      We show (see Theorem~\ref{kummer}) that if an abelian surface
$A$ is $F$-split (hence ordinary), then the associated Kummer $K3$
surface $X$ is also $F$-split (and hence ordinary).  We recall that
the surface $X$ is obtained by blowing up the singularities of the
singular surface $\tilde{A}$, obtained by identifying the points $x$
and $-x$ in $A$.  Although this result can be proved by $l$-adic
methods when the base field is finite, in this paper we prove this for
a perfect field $k$, by relating ordinarity to Frobenius splitting for
such varieties.  In the case of abelian and $K3$-surfaces, it is
possible to compare the notion of ordinarity with that of Frobenius
splitting, and this allows us to handle the passage to the singular
variety $\tilde{A}$.
 
     Ordinary varieties are Hodge-Witt, in that the de Rham-Witt
cohomology groups $H^i(X, W\Omega^j_X)$ are finitely generated over
$W$.  A natural question that arises is whether Frobenius split
varieties are Hodge-Witt. We show (see Theorem~\ref{finiteness}) that
for any smooth projective $F$-split variety over an algebraically
closed field the cohomology groups $H^i(X,W(\O_X))$ are of finite type
as $W$-modules. Using the work of Illusie-Raynaud, we also see that
the first differential $d_1^{i,0}$ is zero for all $i\geq 0$.  During
the course of writing of this paper, the first author refined these
methods to control the nature of crystalline torsion for $F$-split
varieties (see \cite{joshi00b}) and has also shown that any smooth,
projective and Frobenius split threefold is Hodge-Witt.

The foregoing results raise the possibility that Frobenius split
varieties should be Hodge-Witt or even ordinary and indeed the
question of whether or not Frobenius split varieties are ordinary
was raised by V. B. Mehta.  After the first draft of this paper
and \cite{joshi00b} were written and circulated, we found however
that this general expectation, which had been further
strengthened by low dimensional results like
Theorem~\ref{surface-fsord} for surfaces and \cite{joshi00b} (for
threefolds), turns out to be false in higher dimensions. In
\cite{joshi00b} it was shown that any Frobenius split, smooth
projective threefold is Hodge-Witt. It turns out that this is
best possible. We found examples of Frobenius split varieties
which are not ordinary (of dimension at least three) and are not
even Hodge-Witt (dimension at least four).

We also give an example of a smooth fibration $X\to Y$ of smooth,
projective varieties $X$ and $Y$, where the base and fibers are
ordinary, but the total space $X$ is not ordinary. This is in
contrast to the fact that if $X$ is the projective bundle
associated to a vector bundle on $Y$, then $X$ is ordinary if and
only if $Y$ is ordinary. A variant of our method also gives an
example of a variety defined over a number field, whose reduction
modulo all but a finite set of primes is Hodge-Witt (and Frobenius split), 
but which has non-ordinary  reduction for infinitely
many primes.

   One of the other motivating questions of this paper, is a
conjecture of Serre (Conjecture~\ref{serre}) formulated originally in
the context of abelian varieties. Let $K$ be a number field, and let
$X$ be a smooth, projective variety defined over $K$. Then the
conjecture is that there should be a positive density of primes of
$K$, at which $X$ acquires ordinary reduction.  We show that if $X$ is
either an abelian variety or a $K3$ surface defined over $K$,
(Theorem~\ref{positive-density}) then there is a finite extension
$L/K$ of number fields, such that the set of primes of $L$ at which
$X$ has ordinary reduction in the case of $K3$ surfaces, or has
$p$-rank at least two if $X$ is an abelian variety, is of density
one. We note here that a proof of the result for a class of K3
surfaces was also given by Tankeev (see \cite{tankeev95}) under
somewhat restrictive hypotheses.  Our proof mirrors closely the proof
given by Ogus for abelian surfaces. The results of this section are
essentially independent of the contents of the rest of the paper.

 In the final section we study the relationship between ordinariness
and  the torsion in the de Rham-Witt cohomology of varieties and
discuss some questions and conjectures.

\noindent{\bf Acknowledgement:} We would like to acknowledge our debt
of gratitude to V.~B.~Mehta who has explained to us many of his ideas
and insights on Frobenius splitting. This paper was motivated to a
large extent by his questions about crystalline aspects of Frobenius
splitting. We would also like to thank Luc Illusie, Minhyong Kim,
Arthur Ogus, A. J. Parmeshwaran, Douglas Ulmer, Adrian Vasiu for
correspondence, conversations and encouragement.

\section{Preliminaries}\label{preliminaries} 
\subsection{Ordinary varieties}\label{ordinary-subsection}
Let $X$ be a smooth projective variety over a perfect field $k$ of
positive characteristic.  Following Bloch-Kato \cite{bloch86} and
Illusie-Raynaud \cite{illusie83b},  we say that $X$ is ordinary  if
$H^i(X,B^j_X)=0$ for all $i\geq 0,\, j>0$, where  \[ B^j_X= {\rm
image}\left( d:\Omega^{j-1}_X\to \Omega^j_X\right).\] If $X$ is an
abelian variety, then it is known that this definition  coincides with
the usual definition \cite{bloch86}.    By
\cite[Proposition 1.2]{illusie90a}, ordinarity is an  open condition in the
following sense: if $X\to S$ is a smooth, proper family of varieties 
parameterized by $S$, then the set of points $s$ in  $S$, such that the
fiber $X_s$ is ordinary is a Zariski open subset of $S$. Although the
following proposition is well known, we present  it here as an
illustration of the power of this fact.

\begin{proposition} 
    Let $A_{g,1,n}$ denote the moduli space of principally  polarized
abelian varieties $A$ of dimension $g$,  equipped with a level structure
of level $n\geq 5$, and $(n,p)=1$.  Then the set of ordinary points is
open and dense in the moduli of principally polarized abelian varieties
of dimension $g$. 
 \end{proposition} 

\begin{proof}
Since the level $n\geq 5$, it is known that we obtain a fine
moduli space over $\Z\left[\frac{1}{n}\right].$ Further since the
points of the moduli space over $\C$, are uniformized by the
Siegel upper half plane, $A_{g,1,n}$ is an irreducible smooth
variety over $\Z\left[\frac{1}{n}\right].$ Since $p$ is coprime
to $n$, the moduli problem specializes, and we see that
$A_{g,1,n}\otimes_{\Z} {{\rm Spec}}({\bf F}_p)$ is irreducible.
By \cite[Proposition 1.2]{illusie90a}, the ordinary locus is
open. To show it is dense it suffices to prove it is not empty.
But this is easily done by choosing an ordinary elliptic curve
(and there is always one in every characteristic) together with
its principal polarization. Then we can take our ordinary abelian
variety with principal polarization to be the product of this
ordinary elliptic curve and we are done. \end{proof}

\subsection{Cartier Operator} 
   Let $X$ be a smooth proper variety over a perfect field of
characteristic $p>0$, and let $F_X$ (or $F$)  denote the absolute Frobenius of
$X$. We recall a few basic facts about Cartier
operators from \cite{illusie79b}. The first fact we need is that we have
a fundamental exact sequence of locally free sheaves
   \begin{equation}
   0 \to B_X^i \to Z_X^i \xrightarrow{C} \Omega^i_X\to 0,
   \end{equation} where $Z^i_X$ is the sheaf of closed $i$-forms, where
$C$ is the Cartier operator.  The existence of this sequence is the
fundamental theorem of Cartier (see \cite{illusie79b}). 
Since the Cartier operator is also the trace map in Grothendieck
duality theory for the finite flat map $F$, we have  a perfect pairing 
\begin{equation}
F_*(\Omega^j_X)\tensor F_*(\Omega^{n-j}_X) \to \Omega_X^n
\end{equation}
 where $n=\dim(X)$, and the pairing
is given by $(\omega_1,\omega_2)\mapsto C(w_1\wedge w_2)$. This pairing
is perfect and $\O_X$-bilinear (see \cite{mehta87}). 

    In particular, on applying $\Hom(-,\Omega^n_X)$ to the exact
sequence  
\begin{equation} 
    0\to B^n_X \to Z^n_X \to \Omega^n_X \to 0
 \end{equation} we get 
\begin{equation} 
0\to \O_X\to F_*(\O_X) \to B^1_X\to 0 
 \end{equation}

\subsection{$F$-split varieties}\label{f-split-subsection}
    In this section we recall a few basic facts about ordinary and
Frobenius split ($F$-split) varieties. In this section $X$ is a
normal, projective variety over a perfect field $k$. Let $F:X\to X$ be
the Frobenius morphism of $X$.

    Recall that $X$ is {\em $F$-split} if the canonical
exact sequence of sheaves 
\begin{equation}\label{fsplitseq} 
0 \to \O_X \to F_*(\O_X) \to B^1_X \to 0 
\end{equation}
 splits. Note that when $X$ is smooth this is an exact sequence of
locally free $\O_X$-modules.
 
    Frobenius splitting was introduced by Mehta and Ramanathan in
\cite{mehta85b} and a number of remarkable properties were
also investigated in that paper.  It is known for instance that an
abelian variety is $F$-split if and only if it is ordinary in the
usual sense (see \cite{mehta87}).

    In analogy with ordinary varieties, we  now consider the openness of
the  F-split  condition. Let $f:X\to S$ be a smooth projective
morphism. Let $X'$ denote the
fiber product $X\times_{(S,F_S)}S$. Then one has a morphism
   \begin{equation}
   F_{X/S}:X\to X'
   \end{equation}
which is called the relative Frobenius morphism (see \cite{illusie79b}).
The restriction of $F_{X/S}$ to the fibers of $f$ induce the Frobenius
morphism on the fibers.

The following proposition is the relative version of Proposition~9 of
\cite{mehta85b}. 

\begin{proposition}\label{relative-fs}
Let $f:X\to S$ be a smooth projective morphism of schemes in
characteristic $p$. Let $F_{X/S}:X\to X'$
be the relative Frobenius morphism. Then for a point $s\in S$, the
fiber $X_s=f^{-1}(s)$ is $F$-split if and only if the natural map
\begin{equation}
    R^nf_*(\Omega^d_{X/S})\tensor k(s) \to
        R^nf_*\left(F^*_{X/S}(\Omega^d_{X/S})\right) \tensor k(s) 
    \end{equation} is
injective. 
 \end{proposition}

We note that the map 
\begin{equation} F^*_S(R^nf_*(\Omega_{X/S}^d))\to
R^nf_*(F^*_{X/S}(\Omega^d_{X/S})) \end{equation} 
is ${\mathcal O}_S$-linear. 
From this and Proposition~\ref{relative-fs}  we get: 
\begin{proposition} With the notations of Proposition~\ref{relative-fs},
there exists a Zariski open subset $U\subset S$ such that all the
fibers of $f$ over points of $U$ are $F$-split. \end{proposition}

\subsection{Ordinary and $F$-split varieties}
In this subsection we record a key lemma which we need and also record
our proof that Frobenius split surfaces are ordinary.  Our main tool
here is the duality induced by the Cartier operator (see
\cite{mehta87}).

We have the following lemma:

\begin{lemma}\label{b1vanish}
Let $X$ be any smooth projective variety over a perfect field. Assume
$X$ is $F$-split. Then for all $i\geq 0$,
    $$H^i(X,B^1_X)=0.$$
\end{lemma}

\begin{proof}
As $X$ is $F$-split, it follows that $F_*(\O_X)=\O_X\oplus B^1_X$. Hence
$$H^i(X,F_*(\O_X))=H^i(\O_X)\oplus H^i(X,B^1_X).$$
But by the Leray
spectral sequence applied to the Frobenius morphism and the
projection formula we see that 
$$H^i(X,F_*(\O_X))\simeq  H^i(X,\O_X)$$
 and hence
we see that 
$$\dim H^i(\O_X)=\dim H^i(\O_X)+\dim H^i(B^1_X)$$ and so the
lemma is proved.
\end{proof}

\begin{theorem}\label{surface-fsord}
Let $X$ be any smooth projective, $F$-split surface over a perfect
field. Then $X$ is ordinary.
\end{theorem}

\begin{proof}
   By Lemma~\ref{b1vanish} we know that for any $F$-split
variety $X$, $$H^i(X, B^1_X)=0$$ for all $i$. So when $X$ is a
surface we need to check that the same vanishing is also valid for
$B^2_X$. But this is immediate from Serre duality and the
following fact: Cartier operator induces a perfect pairing (we
write it under assumption that $X$ is a surface)
   $$F_*(\O_X)\tensor F_*(\Omega^2_X) \to \Omega^2_X,$$ given by
$(f,w)\mapsto C(fw)$ and this induces a perfect pairing $B^1_X\tensor
B^2_X\to \Omega^2_X=\omega_X$ (see \cite{mehta87}).
\end{proof}

\begin{theorem}\label{ktrivial-ordfs}
Let $X$ be any smooth projective, ordinary variety over a perfect
field $k$. If the canonical bundle of $X$ is trivial, then $X$ is
Frobenius split. 
\end{theorem}

\begin{proof} The obstruction to the splitting of the sequence 
 \begin{equation*}
   0 \to \O_X \to F_*(\O_X) \to B^1_X \to 0,
   \end{equation*}
is an element of ${\rm Ext}^1(B^1_X, \O_X)\simeq  H^1(X, (B^1_X)^*). $ 
The duality pairing induced by the Cartier operator implies that 
$$(B^1_X)^* \simeq  B^n_X\otimes \omega_X, $$
where $\omega_X$ denotes the canonical bundle. 
Since we have assumed that $\omega_X$ is trivial, it follows that
\[{\rm Ext}^1(B^1_X, \O_X)\simeq H^1(X, (B^1_X)^*)\simeq 
H^1(X, B^n_X)=0 \]
where the  vanishing follows from the ordinarity assumption. Hence $X$ 
is Frobenius split.
\end{proof}

\section{De Rham-Witt cohomology of $F$-split varieties}
\label{hodge-witt-section}

\subsection{de Rham-Witt cohomology}
   The standard reference for de Rham-Witt cohomology is
\cite{illusie79b}. Throughout this section,  the
following notations will be in force. Let $k$ be an algebraically
closed field of characteristic $p>0$, and $X$ a smooth, projective
variety over $k$.  Let $W=W(k)$ be the ring of Witt
vectors of $k$. Let $K=W[1/p]$ be the quotient field of $W$. Note that
as $k$ is perfect, $W$ is a Noetherian local ring with a discrete
valuation and with residue field $k$. For any $n\geq 1$, let
$W_n=W(k)/p^n$. $W$ comes equipped with a lift $\sigma:W\to
W$, of the Frobenius morphism of $k$, which will be called the Frobenius
of $W$. We define a non-commutative ring $R^0=W_\sigma[V,F]$, where
$F,V$ are two indeterminate subject to the relations $FV=VF=p$
and $Fa=\sigma(a)F$ and $aV=V\sigma(a)$. The ring $R^0$ is called the
Dieudonne ring of $k$. The notation is borrowed from
\cite{illusie83b}.

     Let $\{W_n\Omega^*_X\}_{n\geq 1}$ be the de Rham-Witt pro-complex
constructed in \cite{illusie79b}. It is standard that for each $n\geq
1,i,j\geq 0$, $H^i(X,W_n\Omega^j_X)$ are of finite type over $W_n$. We
define
\begin{equation}
    H^i(X,W\Omega^j_X)=\liminv_{n}H^i(X,W_n\Omega^j_X),
 \end{equation}  which are $W$-modules of finite type up to torsion.
These cohomology groups are called Hodge-Witt cohomology groups  of $X$.

\begin{definition} $X$ is Hodge-Witt if for $i,j\geq 0$, 
the Hodge-Witt cohomology
 groups $ H^i(X,W\Omega^j_X)$ are finite type $W$-modules.

\end{definition}

    The properties of the de Rham-Witt pro-complex are reflected  in
these cohomology modules and in particular we note that for each
$i,j$, the Hodge-Witt groups $H^i(X,W\Omega_X^j)$ are left modules
over $R^0$. The complex $W\Omega^*_X$ defined in a natural manner from
the de Rham-Witt pro-complex computes the crystalline cohomology of
$X$ and in particular there is a spectral sequence
\begin{equation}
   E_1^{i,j}=H^i(X, W\Omega^j_X)\Rightarrow H^*_{\rm cris}(X/W)
\end{equation}

    This spectral sequence induces a filtration on the crystalline
cohomology of $X$ which is called the slope filtration and the
spectral sequence above is called the slope spectral sequence (see
\cite{illusie79b}). It is standard (see \cite{illusie79b} and
\cite{illusie83b}) that the slope spectral sequence degenerates
at $E_1$ modulo torsion (i.e. the differentials are zero on tensoring
with $K$) and at $E_2$ up to finite length (i.e. all the differential
have images which are of finite length over W). 

   In dealing with the slope spectral sequence it is more
convenient to work with a bigger ring than $R^0$. This ring was
introduced in \cite{illusie83b}.  Let $R=R^0\oplus R^1$ be a
graded $W$-algebra which is generated in degree $0$ by variables $F,V$
with the properties listed earlier (so $R^0$ is the the Dieudonne ring
of $k$) and $R^1$ is a bimodule over $R^0$ generated in degree $1$ by
$d$ with the properties $d^2=0$ and $FdV=d$, and $da=ad$ for any $a\in
W$. The algebra $R$ is called the Raynaud-Dieudonne ring of $k$ (see
\cite{illusie83b}). The complex $(E_1^{*,i},d_1)$ is a graded
module over $R$ and is in fact a coherent, left $R$-module (in a
suitable sense, see \cite{illusie83b}).

\subsection{A finiteness result}
    For a general variety $X$, the de Rham-Witt cohomology groups are
not of finite type over $W$, and the structure of these groups
reflects the arithmetical properties of $X$. For instance, in
\cite{bloch86}, \cite{illusie83b} it is shown that for ordinary
varieties $H^i(X,W\Omega^j)$ are of finite type over $W$. The
following theorem lends more evidence towards the general expectation
that Frobenius split varieties should be ordinary.

\begin{theorem}\label{finiteness}
   Let $X/k$ be any smooth, projective, $F$-split variety over an
algebraically closed field $k$ of characteristic $p>0$. Then for each
$i\geq 0$, $H^i(X,W(\O_X))$ is a finite type $W(k)$-module.
\end{theorem}

\begin{remark}
      When the formal Brauer group $\widehat{{Br}^i(X)}$ of $X$
(associated by Artin and Mazur; see \cite{artin77}) is representable,
$H^i(X,W(\O_X))$ is the Cartier module of this formal group. When this
module is free of finite type over $W$, the formal Brauer group is a
$p$-divisible group of height equal to the dimension 
$\dim_K H^i(X,W(\O_X))\tensor_WK$.  \end{remark}

\begin{corollary}
   Let $X$ be as in Theorem~\ref{finiteness}. Then for all
$i\geq0$, the differential
\begin{equation}
d_1^{i,0}: H^i(X,W(\O_X)) \to H^i(X,W\Omega_X^1)
\end{equation}
is zero.
\end{corollary}

\begin{proof}
   This is immediate from Theorem~\ref{finiteness} and
\cite{illusie83b}.
\end{proof}

\begin{remark}
   Let $X$ be any smooth projective variety. In
\cite{illusie83b} it is has been shown that for all $j$,
$H^1(X,W\Omega^j_X)$ are finite type $W$-modules.
\end{remark}

\begin{remark} By combining the theory of dominoes with 
Theorem~\ref{finiteness} and the foregoing remark, it can be shown
that Frobenius split, smooth, projective threefolds are Hodge-Witt
\cite{joshi00b}.
\end{remark}

      Before we give the proof of the above theorem, we need a few
preparatory lemmas.  Our proofs use the theory of higher Cartier
operators as outlined in \cite{illusie79b}. To set up conformity with
notations from previous sections we recall Illusie's notations and our
definitions
\begin{equation} B_1\Omega^1_X=B^1_X={\rm
image}(d:\O_X=\Omega^0_X\to \Omega^1_X) \end{equation}
and
\begin{equation}
Z_1\Omega^1_X=Z^1={\rm ker}(d:\Omega^1_X\to \Omega^2_X)
\end{equation}

The higher Cartier sheaves, $Z_n\Omega^i_X,B_n\Omega^i_X$ are defined
inductively in \cite{illusie79b} (see page 519 of \cite{illusie79b}). The
formation of these sheaves is compatible with arbitrary base change
(see \cite{illusie79b}, page~519).

\begin{lemma} \label{vanishingbn1}
Let $X$ be any smooth, projective, $F$-split variety over
$k$. Then for all $n\geq 0$ and for all $i\geq 0$ we have
$H^i(X,B_n\Omega^1_X)=0$. \end{lemma} 

\begin{proof} 
   The case $n=1$ is Lemma~\ref{b1vanish}.  We prove the result
by induction on $n$. We recall that we have an exact sequence (the
arrow on the extreme right is the Cartier operator):
   $$ 0 \to B_1\Omega^1_X \to B_{n+1}\Omega^1_X \xrightarrow{C^{-1}}
B_n\Omega^1_X\to 0.$$ (see \cite{illusie79b}, page 519). This is
essentially the definition of $B_{n+1}$ using $B_n$. The result now
follows trivially from the above exact sequence and the result for
$n=1$. \end{proof}

\begin{lemma}
    Let $X$ be a smooth, projective, $F$-split variety over an
algebraically closed field $k$ of characteristic $p>0$. Then for all
$n\geq 0$ and for all $i\geq 0$, we have $H^i(X,Z_n\Omega^1_X)\simeq
H^i(X,\Omega^1_X)$. In particular we have
$\forall i>0,\forall n\geq 1$,
$$\dim H^i(X,Z_n\Omega^1_X)= \dim H^i(X,\Omega^1_X).$$
\end{lemma}

\begin{proof}
We recall the exact sequence (page 531, 2.5.1.2 of \cite{illusie79b}):
   $$0\to B_n\Omega^1_X \to Z_n\Omega^1_X \to \Omega^1_X \to 0.$$
Then proof follows from the vanishing of cohomology of $ B_n\Omega^1_X$.
\end{proof}

\begin{proof}{[of Theorem~\ref{finiteness}]} 
      By (\cite{illusie79b}, page~613, Proposition~2.16), it
suffices to prove that for all $j\geq 0$ and for all $n\geq 0$,
$H^j(X,Z_n\O_X)$ has bounded dimension. But
$Z_1\O_X=\ker(d:\O_X\to\Omega^1_X)$ and by definition $Z_n\O_X\to
Z_{n-1}\O_X$ is an isomorphism for all $n\geq 2$ (the arrow in this
isomorphism is the Cartier operator, (\cite{illusie79b}, see 2.5.1.2,
page 531). Thus the required cohomology has dimension independent of
$n$. Next we need to check $H^j(X,B_n\Omega^1_X)$ has bounded
dimension for all $n$. But by Lemma~\ref{vanishingbn1}, this group is
zero! Thus we can apply Proposition 2.16 of \cite{illusie79b} to
deduce that $H^i(X,W(\O_X))$ is a finite type $W(k)$-module.
\end{proof}

\section{Examples}
\subsection{Two Questions}
One of the main motivations for this paper, are the following
questions  raised by
V.~B.~Mehta.
\begin{question}
  Is any smooth projective, Frobenius split variety over a perfect
field of characteristic $p$ of Hodge-Witt type?
\end{question}

The above question is a weaker variant of the following.

\begin{question}
   Is any smooth projective, Frobenius split variety Bloch-Kato
ordinary?
\end{question}

    We know by Theorem \ref{surface-fsord} that a smooth,
projective Frobenius split surface is ordinary. In \cite{joshi00b}, it
is shown that any Frobenius split smooth,
projective three fold is Hodge-Witt. Further it is
known that for abelian varieties, the notions of Frobenius splitting
and ordinarity coincide \cite{mehta87}.  However in contrast to the
expectation created by these
results, we show in this section that the first question is false in
dimensions greater than 3, and  the second question is false in 
dimensions bigger than two.

   Our examples also give examples of varieties which are
Hodge-Witt, but are not ordinary. These examples also provide examples of
smooth, projective varieties $f:~ X\to Y$, such that both $Y$ and the
the (smooth) fibers of $f$ are ordinary, but $X$ is not ordinary.

\subsection{}\label{frobenius-splitting-section}
   Let $X$ be a smooth, projective variety with canonical bundle
$\omega$.  We recall that by Cartier duality there is functorial
isomorphism \cite[Proposition 5]{mehta85b}, \cite{mehta91},
\[ F_*\omega^{1-p}\simeq {\rm Hom}_{\O_X}(F_*\O_X, \O_X).\]
By means of this, we obtain a natural identification
\[ H^0(X, \omega_X^{1-p})\simeq {\rm Hom}_{\O_X}(F_*\O_X, \O_X).\]

\begin{definition} A section $\sigma \in H^0(X, \omega_X^{1-p})$, such that
under the above isomorphism, $\sigma$ provides a splitting of $X$, will be
called  as a splitting section. 
\end{definition}

      The key result we need is a criteria on the relative embedding
of a smooth subvariety in a smooth, Frobenius split variety, such that
the blow up along the subvariety remains Frobenius split.  We recall
now some of the concepts and results regarding compatible Frobenius
splitting of subvarieties.  Let $X$ be a Frobenius split variety, and
let
\[ \sigma: F_*({\mathcal O}_X)\to {\mathcal O}_X\]
be a splitting of the Frobenius morphism. Suppose $Y$ is a subvariety 
of $X$, defined by a sheaf of ideals ${\mathcal I}_Y\subset 
{\mathcal O}_X$. In this case we have a notion of $Y$ being compatibly 
Frobenius split in $X$ as follows:

\begin{definition} $Y$ is said to be compatibly split by $\sigma$ in $X$, if 
$$\sigma(F_*({\mathcal I}_Y))\subset {\mathcal I}_Y.$$
\end{definition}

   Let $X$ be a nonsingular variety and $Y$ a nonsingular
subvariety of codimension $d\geq 2$.  Denote by $B_Y(X)$ the blow up
of $X$ along $Y$, and by $E$ the exceptional divisor.  The following
result follows quite easily from \cite[Proposition
2.1]{lakshmibai98}.

\begin{proposition}
Let  $s\in H^0(X, \omega_X^{-1})$. Suppose that $s^{p-1}$ is a
splitting section of $X$, and that it vanishes to order   
$(d-1)$ or $d$ generically along $Y$. 
Then $s^{p-1}$ extends to  a splitting of $B_Y(X)$.
Moreover if $s$ vanishes to order $d$ generically along $Y$, then $E$
is compatibly split in $B_Y(X)$.
\end{proposition}
 
   On the other hand, we have the following criterion for a
blow-up to be ordinary or Hodge-Witt \cite{illusie90a}, \cite{gros85}:
\begin{proposition} $B_Y(X)$ is ordinary (or Hodge-Witt) if and only if
both $X$ and $Y$ are ordinary (resp. Hodge-Witt).
\end{proposition} 
The proof of this proposition follows from the decomposition of
$W$-modules, compatible with the action of the Frobenius \cite[IV
1.1.9]{gros85},
\[ H^j(X, W\Omega^i_X)\oplus(\oplus_{0<l<d}H^{j-l}(Y,W\Omega_Y^{i-l}))
\tilde{\to} H^j(B_Y(X), W\Omega_{B_Y(X)}^i), \] 
and the fact that a smooth, proper variety $Z$ is ordinary if and only if 
\[F: H^j(Z, W\Omega_Z^i)\to H^j(Z, W\Omega_Z^i)\]
is an isomorphism for all $i$ and $j$. 

Combining the above propositions, we obtain the following theorem: 

\begin{theorem}\label{blowupthm} 
   Let $X$ be a smooth, projective Frobenius split variety, with
a splitting section $s^{p-1}$ as above.  Suppose that $s$ vanishes to
order precisely $(d-1)$ generically along a smooth subvariety $Y$ of
codimension $d$ in $X$. Further assume that $Y$ is not ordinary (or
not Hodge-Witt). Then $B_Y(X)$ is Frobenius split but not ordinary
(resp. not Hodge-Witt).
\end{theorem} 

\subsection{} 
We can now give the examples of Frobenius split varieties which are
not ordinary or Hodge-Witt.

\begin{example}\label{ellipticase} 
   Let $E$ be a supersingular elliptic curve in ${\mathbf
P}^3$. $E$ is contained in the zero locus of  non-degenerate quadric
$q$. Let $t_1$ and $t_2$ be linear polynomials such that the ideals
generated by choosing any combination of $q,~t_1, ~t_2$ define
complete intersection subvarieties in ${\mathbf P}^3$. Then it can be
checked using \cite[Proposition 7]{mehta85b} that the section
$s=qt_1t_2\in \O(4)$ gives rise to a splitting of ${\mathbf P}^3$,
vanishing to order $1$ along $E$. Hence the blow up of ${\mathbf P}^3$
is Frobenius split (and is Hodge-Witt) but is not ordinary.
\end{example}

\begin{example}\label{pencil} 
   A natural question that arises in the study of the geometry of
ordinary varieties, is whether a variety is ordinary, if it is 
fibered over an ordinary
variety, such that  the smooth fibres are ordinary.  
The above example also provides an example of a variety fibered
over ${\mathbf P}^1$, such that the (smooth) fibers are ordinary but the
variety itself is not ordinary. The elliptic curve is defined as the
complete intersection of two non-degenerate quadrics which generates a
pencil of quadrics. The strict transform of these quadrics in the
blow-up gives a fibration of $B_E({\mathbf P}^3)$ over ${\mathbf P}^1$
by Frobenius split varieties. However Frobenius split (smooth) surfaces are
ordinary, and this gives us the desired example. 
\end{example}

\begin{example}\label{highdim}
    The above example can be generalized. Let $Y$ be a smooth
hypersurface in ${\mathbf P}^{n+1}\subset {\mathbf P}^{n+2}$, for
example a Fermat hypersurface of degree $m$. Choose a system of
coordinates $x_0,\cdots, x_{n+1}$ on ${\mathbf P}^{n+1}$, where
${\mathbf P}^{n+1}$ is given by $x_0=0$. Then $s^{p-1}=(x_0\cdots
x_{n+2})^{p-1}$ is a splitting section vanishing precisely to order
$(p-1)$ generically along $Y$. The blow-up of ${\mathbf P}^{n+2}$
along $Y$ is then Frobenius split. Recall that results of
\cite{toki96} give explicit conditions on $(n,p,m)$ under which the
Fermat hypersurface $Y$ is not Hodge-Witt. For instance assume that
$p\not \equiv 1 \mod{m}$, $p$ does not divide $m$ and $n\geq
5,m\geq5$. Then this hypersurface is not Hodge-Witt and so the blowup
is not Hodge-Witt. But the blowup is Frobenius split but neither
Hodge-Witt nor ordinary. The results of \cite{toki96} can also be used
to give examples in dimensions four and five as well.
\end{example} 

\begin{example}\label{numberfieldcase}
   Let $E/\Q$ be an elliptic curve. We embed $E$ in $\P^3$ by
using the embedding given by the linear system $4(\infty)$. The blowup
of $\P^3$ along $E$ has $F$-split and Hodge-Witt reduction at all but
finite number of primes. In fact, by the blowup formula for Hodge-Witt
cohomology, blowup of $\P^3$ along any smooth projective embedded
curve is Hodge-Witt. By \cite{elkies87} we know that the reduction of
$E$ is supersingular at infinitely many primes, and if $E$ has CM,
then it has supersingular reduction at a set of primes of density
$1/2$.   Hence there are 
infinitely many primes where the blowup has $F$-split (and Hodge-Witt)
but non-ordinary reduction.

\end{example}

\begin{example} The above examples might lead one to raise the
question whether  Fano $F$-split varieties are ordinary or Hodge-Witt.
But even this turns out to be false.  See Example \ref{fano}.
\end{example}

\section{Kummer surfaces}\label{kummer-section}
In this section we explore further the relationship between Frobenius split
and ordinary varieties, especially in the context of Kummer
surfaces.

\subsection{Kummer Surfaces over perfect fields}
Let $A$ be an abelian surface over a perfect field $k$ of odd,
 positive characteristic $p$. Denote by $\iota:A\ra A$, the involution
 $x\mapsto -x$ on the abelian surface. Let $\tilde{A}$ denote the
 quotient variety of $A$ with respect to this involution. It is known
 that $\tilde{A}$ is Gorenstein having only quotient
 singularities. Denote by $\omega_{\tilde{A}}$ the dualizing sheaf. It
 is known that $\omega_{\tilde{A}}\simeq {\mathcal O}_{\tilde{A}},$
 the structure sheaf on $\tilde{A}$.

There exists a smooth, projective variety $X$, which is a blow up of 
$\tilde{A}$, at the sixteen singular points of $\tilde{A}.$ $X$ is 
a $K3$-surface, in that it is simply connected and $H^2(X, {\mathcal O}_X)$
is one dimensional. $X$ is the Kummer $K3$-surface associated to the abelian 
surface $A.$

\begin{theorem}\label{kummer-fs} With notation as in the above theorem,  
$A$ is Frobenius split if and only if $X$ is Frobenius split.
\end{theorem}

 Since the canonical bundles of
$A$ and $X$ are trivial, we see by Theorems \ref{surface-fsord} and
\ref{ktrivial-ordfs}, that the above theorem is equivalent to proving
the following:

\begin{theorem}\label{kummer}
Let $A$ be an abelian surface over a perfect field $k$ of characteristic 
$p>0$. Let $X$ be the associated Kummer $K3$ surface. 
Then $X$ is ordinary if and only if $A$ is ordinary. 
\end{theorem}

\begin{proof}
We will need the formalism of Section~\ref{frobenius-splitting-section} for
Gorenstein varieties. We recall this from \cite{mehta91}.

   Let $Z$ be a scheme over $k$. In order to split $Z$, we need a
map $F_*\O_Z\to \O_Z$ such that the composite with the $\O_Z\to
F_*\O_Z\to \O_Z$ is the identity. Suppose now that $Z$ is a reduced
equidimensional Gorenstein $k$-scheme. By applying duality for the
Frobenius morphism, we obtain a canonical isomorphism of sheaves on
$Z$, as in \cite[Lemma 1]{mehta91},
\begin{equation}\label{F-duality} 
F_*\omega_Z^{1-p}\simeq {\rm Hom}_{\O_Z}(F_*\O_Z, \O_Z),
\end{equation}
where $\omega_Z$ denotes the dualizing sheaf of $Z$.  In particular,
Frobenius splittings of $Z$ are induced by sections of $H^0(Z,
\omega^{1-p}_Z)$.

Let $Z$ denote any one of the varieties $A,\, \tilde{A},\, X$. By
definition $A$ and $X$ are smooth, and it is known that $\tilde{A}$ is
Gorenstein. Further the dualizing sheaf of $Z$ is the structure sheaf
$\O_Z$ in each of the above cases. Moreover the $\Z/2\Z$ action
is trivial on $H^0(A,\O_A)$, we have a natural isomorphism
\[H^0(A,\O_A)\simeq H^0(\tilde{A}, \O_{\tilde{A}})\simeq H^0(X,\O_X).\]
Let $s$ denote a section in any one of the above cohomology groups, and 
we continue to denote by $s$, its image in the other cohomology groups. 
By the isomorphism \ref{F-duality}, $s$ gives rise to a morphism  
$F_*\O_Z\to \O_Z.$ To check that $s$ gives  a splitting section,
that the composite $\O_Z\to F_*\O_Z\to \O_Z$ is the identity, it is 
enough to check at a point $P$ on $Z$, since $Z$ is projective and any 
global map $\O_Z\to \O_Z$ is a constant. By the local nature of duality,
the morphism \ref{F-duality}, is an isomorphism of sheaves, and it is 
enough to check the  splitting condition in the formal neighborhood of 
a smooth point $P$ on $Z$.
 
   We now choose  $P$ to be a non $2$-torsion point on $A$. We
continue to denote by $P$, the image of $P$ in $\tilde{A}$ and $X$. We
then have an isomorphism of the formal completions,
$$\hat{\O}_{A, P}\simeq \hat{\O}_{\tilde{A}, P}\simeq \hat{\O}_{X, P}$$
compatible with the isomorphism \ref{F-duality}. Hence a section $s$ 
gives a splitting section for $A$ if and only if it gives a splitting 
section for $\tilde{A}$,  or equivalently  for $X$. 
Hence $A$ is $F$-split is equivalent 
to $X$ being  $F$-split, and this is equivalent to   
$\tilde{A}$ being $F$-split.
\end{proof} 

\begin{remark} Over finite fields, it is possible to give a different
proof using $l$-adic methods.  For a smooth, projective surface $X$ with trivial canonical 
bundle, the condition for being Frobenius split is that the Frobenius
$F:H^2(X, \O_X)\to H^2(X, \O_X)$ is an isomorphism. When the surfaces
are defined over finite fields, it follows from the Katz congruence
formula for the zeta function \cite{katz73c}, applied to a
$K3$-surface, that there is precisely one eigenvalue of the
crystalline Frobenius acting on $H^2_{cris}(X/W)$ which is a $p$-adic
unit. From the shape of the Hodge polygon and duality in the case of
abelian and $K3$-surfaces, we then conclude that this is equivalent to
$K3$-surface being ordinary. By comparing $H^2(A,\O_A)$ and
$H^2(X,\O_X)$, we obtain a different proof over finite
fields, that the Kummer $K3$-surface is ordinary if the abelian
surface is ordinary. However these $l$-adic methods do not seem to
generalize to an arbitrary perfect base field.  
\end{remark}

\begin{remark} In the course of the proof of the Tate conjecture for
ordinary $K3$-surfaces $X$ over a finite field \cite{nygaard83}, it is
shown that the
Kuga-Satake abelian variety $K(X)$ associated to $X$ is ordinary,
provided  $X$ is
ordinary. If $X$ is a Kummer $K3$-surface associated to an abelian
surface $A$, then it is known that $K(X)$ is isogeneous to a sum of
copies of $A$. It follows that if $X$ is ordinary, then $A$ is
ordinary.  
\end{remark} 

\begin{remark} Combining the above theorem with a theorem of Ogus
\cite[page 372]{ogus82}, it follows that for a Kummer $K3$-surface 
defined over a number field, there is a finite extension over which the 
variety acquires ordinary reduction at a set of primes of density one. 
In the next section, we will show that Ogus' proof 
extends to prove this result for any $K3$-surface defined over a number 
field.
\end{remark}

\section{Primes of ordinary reduction for $K3$ surfaces}   
   The following more general question, which is one of the
motivating questions for this paper, is the following conjecture which
is well-known and was raised initially for abelian varieties by Serre:
\begin{conj}\label{serre}
   Let $X/K$ be a smooth projective variety over a number field. Then
there is a positive density of primes $v$ of $K$ for which $X$ has
good ordinary reduction at $v$.
\end{conj}

Let $K$ be a number field, and let $X$ denote either an abelian
variety or a $K3$ surface defined over $K$. Our aim in this section is
to show that there is a finite extension $L/K$ of number fields, such
that the set of primes of $L$ at which $X$ has ordinary reduction in
the case of $K3$ surfaces, or has $p$-rank at least two if $X$ is an
abelian variety of dimension at least two, is of density one. Our
proof closely follows the method of Ogus for abelian surfaces (see
\cite[page 372]{ogus82}).

   We note here that a proof of the result for a class of K3
surfaces was also given by Tankeev (see \cite{tankeev95}) under some
what restrictive hypothesis. The question of primes of ordinary
reduction for abelian varieties has also been treated recently by
R.~Noot (see \cite{noot95}), R.~Pink (see \cite{pink98}) and more
recently A.~Vasiu (see \cite{vasiu00}) has studied the question for a
wider class of varieties. The approach adopted by these authors is
through the study of Mumford-Tate groups.

   Let $\O_K$ be the ring of integers of $K$; for a finite place
$v$ of $K$ lying above a rational prime $p$, let $\O_{v}$ be the
completion of $\O_K$ with respect to $v$ and let $k_v$ be the residue
field at $v$ of cardinality $q_v=p^{c_v}$.  Assume that $v$ is a place
of good reduction for $X$ as above and write $X_v$ for the reduction
of $X$ at $v$.  We recall here the following facts:

\subsection{} (Weil, Deligne,  Ogus)\cite{ogus82}: 
   The Frobenius endomorphism $F_v$ is a semi-simple endomorphism
 of the $l$-adic cohomology groups $H^i_l:= H^i_{et}(X\otimes
 \overline{K}, \Q_l)$ for a prime $l\neq p$.  The $l$-adic
 characteristic polynomial $P_{i, v}(t)= {\rm det}(1-tF_v\mid
 H^i_{et}(X\otimes \overline{K}, \Q_l)$ is an integral polynomial and
 is independent of $l$. Let
\[ a_v={\rm Tr}(F_v\mid H^2_{et}(X\otimes \overline{K}, \Q_l)\]
denote the trace of the $l$-adic Frobenius acting on the second \'etale 
cohomology group. $a_v$ is a rational integer.  
 
\subsection{} (Deligne-Weil estimates) \cite{deligne74a}: 
It follows from Weil estimates proved by Weil for abelian varieties and
by Deligne in general  that
\[ |a_v|\leq dp\]
where $d=\dim H^2_l$ is a constant independent of the place $v$. 

\subsection{} (Katz-Messing theorem) \cite{katz74}:
 Let $\phi_v$ denote the 
crystalline Frobenius on $H^i_{cris}(X/W(k_v))$. $\phi_v^{c_v}$ is
linear over $W(k_v)$, and  the characteristic polynomials of 
the crystalline Frobenius and
the $l$-adic Frobenius are equal:  
$$P_{i,v}(t)={\rm det}(1-t\phi_v^{c_v}\mid 
H^i_{cris}(X/W(k_v))\otimes K_v).$$  

\subsection{}(Semi-simplicity of the crystalline
Frobenius)\cite{ogus82}: If $X$ is a $K3$-surface, then it is known by
\cite[Theorem 7.5]{ogus82} that the crystalline Frobenius
$\phi_v^{c_v}$ is semi-simple. We remark that although this result is
not essential in the proof of the theorem, it simplifies the proof.

\subsection{}(Mazur's theorem) \cite{mazur72}, \cite{deligne81a}, 
\cite{berthelot-ogus}: 
   There are two parts to the theorem of Mazur that we require. 
After inverting finitely many primes $v\in S$ in $K$, we can assume
that $X$ has good reduction outside $S$. Using Proposition~\ref{fintor}
(see below) we can assume that
$H^i(X, \Omega^j_X)$ and $H^i(X/W(k_v))$ 
are torsion-free outside a finite set of primes of $K$. 
As $X$ is defined over characteristic
zero, the Hodge to de Rham spectral sequence degenerates at $E_1$
stage. Thus all the hypothesis of Mazur's theorem are satisfied. The
two parts of Mazur's theorem that we require are the following:

\subsubsection{}\label{mazur}(Mazur's proof of  Katz's conjecture): 
    Let $L_v$ be a finite extension of field of fractions of
$W(k_v)$, over which the polynomial $P_{i,v}(t)$ splits into linear
factors. Let $w$ denote a valuation on $L_v, $ such that
$w(p)=1$. The Newton polygon of the polynomial $P_{i,v}(t)$ lies above
the Hodge polygon in degree $i$, 
 defined by the Hodge numbers $h^{j, i-j}$ of degree
$i$. Moreover they have the same endpoints. 

\subsubsection{} \label{divisibility}
(Divisiblity property): The crystalline Frobenius $\phi_v$ is divisible 
by $p^i$ when restricted to 
$F^iH^j_{dR}(X/W(k_v)):= \H(X/W(k_v)), \Omega^{\geq i}_X)$.

\subsection{}(Crystalline torsion):
We will also need the following proposition which is certainly
well-known but as we use it in the sequel, we record it here for
convenience.
\begin{proposition}\label{fintor}
Let $X/K$ be a smooth projective variety. Then for all but finitely
many nonarchimedean places $v$, the crystalline cohomology $H^i_{\rm
cris}(X_v/W(k_v))$ is torsion free for all $i$.
\end{proposition}

\begin{proof} 
   We choose a smooth model $\mathcal{X}\to \Spec(\O_K)-V(I)$ for
some non-zero proper ideal $I\subset \O_K$. The relative de Rham
cohomology of the smooth model $\mathcal{X}$ is a finitely generated
$\O_K$-module and has bounded torsion.  After inverting a finite set
$S$ of primes, we can assume that $H^i_{dR}({\mathcal X}, \O_S)$ is a
torsion-free ${\mathcal O}_S$ module, where ${\mathcal O}_S$ is the
ring of $S$-integers in $K$.  By the comparison theorem of Berthelot
(see \cite{berthelot-cohomologie}), there is a natural isomorphism of
the crystalline cohomology of $X_p$ to that of the de Rham cohomology
of the generic fiber of a lifting to $\Z_p$.
$$ H^i_{\rm cris}(X_v/W(k_v)) \equiv H^i_{dR}({\mathcal X},
\O_S)\otimes W(k_v).$$
This proves our proposition. Moreover the proof
shows that we can assume after inverting some more primes, that the
Hodge filtrations $F^jH^i_{\rm dR}({\mathcal X}, \O_S)$ are also
locally free over $\O_S$, such that the sub-quotients are also locally
free modules.
\end{proof}

We first note the following lemma which is fundamental to the proof.

\begin{lemma} With notation as above, assume the following:\\
a) if $X$ is a $K3$-surface, then $X$ does not have ordinary reduction
at $v$. \\
 b) if $X$ is an abelian variety, then the $p$-rank of the
reduction of $X$ at $v$, is at most 1.
Then $p|a_v.$
\end{lemma}

\begin{proof} 
   Let $w$ be a valuation as in \ref{mazur} above. If $X$ is an
abelian variety defined over the finite field $k_v$, then the $p$-rank
of $X$ is precisely the number of eigenvalues of the correct power of
the crystalline Frobenius acting on $H^1(X/W(k_v))\otimes \Q_p$, which
are $p$-adic units.  Suppose now $\alpha$ is an eigenvalue of the
crystalline Frobenius $\phi_v^{c_v}$ acting on $H^2(X/W(k_v))\otimes
\Q_p.$ In case b), the hypothesis implies that $w(\alpha)$ is
positive, and hence $w(a_v)$ is strictly positive. As $a_v$ is a
rational integer, the lemma follows.

   When $X$ is a $K3$-surface, it follows from the shapes of the
Newton and Hodge polygons, that ordinarity is equivalent to the fact
that precisely one eigenvalue of $\phi_v^{a_v}$ acting on
$H^2(X/W(k_v))\otimes \Q_p$ is a $p$-adic unit. Hence if $X$ is not
ordinary, then for any $\alpha$ as above, we have $w(\alpha_v)$ is
positive. Again since $a_v$ is a rational integer, the lemma follows.
\end{proof}

\begin{theorem}\label{positive-density}
   Let $X$ be a $K3$ surface or an abelian variety of dimension at
least two defined over a number field $K$. Then there is a finite
extension $L/K$ of number fields, such that
\begin{enumerate}
   \item if $X$ is a $K3$-surface, then $X\times_K L$ has
      ordinary reduction at a set of primes of density one in $L$.
   \item if $X$ is an abelian variety, then there is set of
      primes $O$ of density one in $L$, such that the reduction of
   $X\times_{K}L$ at a prime $p\in O$ has $p$-rank at least two.
\end{enumerate}
 \end{theorem}

\begin{proof} 
   Our proof follows closely the method of Serre and Ogus (see
\cite{ogus82}).  Fix a prime $l$, and let $\rho_l$ denote the
corresponding Galois representation on $H^2_l$. The Galois group $G_K$
leaves a lattice $V_l$ fixed, and let $\overline{\rho}_l$ denote the
representation of $G_L$ on $V_l\otimes \Z/l\Z.$ Let $L$ be a Galois
extension of $\Q, $ containing $K$ and the $l^{th}$ roots of unity,
and such that for $\sigma \in G_L$, $\overline{\rho}_l(\sigma)=1.$ We
have
   $$a_v\equiv d ({\rm mod}~l),$$
where $d={\rm dim}~H^2_l$. 

    Let $v$ be a prime of $L$ of degree $1$ over $\Q,$ lying over
the rational prime $p.$ Since $p$ splits completely in $L$, and $L$
contains the $l^{th}$ roots of unity, we have $p\equiv 1 ({\rm
mod}~l).$ Now choose $l>d$. Since $|a_v|\leq dp$ and is a rational
integer divisible by $p$ from the above lemma, it follows on taking
congruences modulo $l$ that
   $$ a_v=\pm dp.$$
Now $a_v$ is the sum of $d$ algebraic integers each of which is 
of absolute value $p$ with respect to any embedding. It follows 
that all these eigenvalues must be equal, and equals $\pm p.$ 
Hence we have that 
   $$F_v=\pm pI$$ 
as an operator on $H^2_l$. By the semi-simplicity of the crystalline
Frobenius for abelian varieties and $K3$-surfaces \cite{ogus82}, it
follows that $\phi_v=\pm pI$. But this
contradicts the divisiblity property of the crystalline 
Frobenius \ref{divisibility}, that the crystalline 
Frobenius is divisible by $p^2$ on  $F^2H^2_{dR}(X\times k_v)$. Hence 
$v$ has to be a prime of ordinary reduction, and this completes the
proof of our theorem.
\end{proof}

   Ogus' method can in fact be axiomatized to give positive
density results whenever certain cohomological conditions are
satisfied. We present this formulation for the sake of completeness.

\begin{proposition}
   Let $X$ be a smooth projective variety over a number field
$K$. Assume the following conditions are satisfied:
\begin{enumerate}
\item $\dim H^2(X,\O_X)=1$,
\item The action of the crystalline Frobenius of the reduction $X_\wp$
of $X$ at a prime $\wp$ is semi-simple for all
but finite number of primes $\wp$ of $K$.
\end{enumerate}
   Then the Galois representation $H^2_{et}(X,\Q_\ell)$ is
ordinary at a set of primes of positive density in $K$ and the
$F$-crystal $H^2_{cris}(X_\wp/W)$ is ordinary for these primes. In
other words, the motive $H^2(X)$ has ordinary reduction for a positive
density of primes of $K$.
\end{proposition}

\section{Primes of Hodge-Witt reduction}
\subsection{Hodge-Witt reduction}
   Let $X$ be a smooth projective variety over a number field $K$. We
fix a model $\mathcal{X} \to \Spec(\O_K)$ which is regular, proper and
flat and which is smooth over a suitable non-empty subset of
$\Spec(\O_K)$. All our results are independent of the choice of the
model. In what follows we will be interested in the smooth fibers of
the map $\mathcal{X} \to \Spec(\O_K)$, in other words we will always
consider primes of good reduction. Henceforth $p$ will always denote
such a prime and the fiber over this prime will be denoted by $X_{p}$.

Since ordinary varieties are Hodge-Witt, we can formulate a weaker
version of Conjecture~\ref{serre}.

\begin{conj}\label{positive-density-hodge-witt}
   Let $X/K$ be a smooth projective variety over a number field
then $X$ has Hodge-Witt reduction modulo a set of primes of $K$ of
positive density.
\end{conj}

    For surfaces the geometric genus appears to detect the size
of the set of primes which is predicted in
Conjecture~\ref{positive-density-hodge-witt}.

\begin{theorem}\label{pgtheorem}
   Let $X$ be a smooth, projective surface with $p_g(X)=0$,
defined over a number field $K$.  Then for all but finitely many
primes $\wp$, $X$ has Hodge-Witt reduction at $\wp$.
\end{theorem}

\begin{proof}
By the results of \cite{nygaard79b}, \cite{illusie79b},
\cite{illusie83b}, it suffices to verify that
$H^2(X_{\wp},W(\O_{X_\wp}))$ is zero for all but finite number of
primes $\wp$ of $K$. But the assumption that $p_g(X)=0$ entails that
$H^2(X,\O_X)=H^0(X,K_X)=0$. Hence by the semicontinuity theorem, for
all but finite number of primes $\wp$ of $K$, the reduction $X_\wp$
also has $p_g(X_\wp)=0$. Then by \cite{illusie79b}, \cite{nygaard79b}
one sees that $H^2(X_\wp,W(\O_{X_\wp}))=0$ and so $X_\wp$ is
Hodge-Witt at any such prime.
\end{proof}

\begin{corollary}
Let $X/K$ be an Enriques surface over a number field $K$. Then $X$ has
Hodge-Witt reduction modulo all but finite number of primes of $K$.
\end{corollary}

\begin{proof}
This is immediate from the fact that for an Enriques surface over $K$,
$p_g(X)=0$. 
\end{proof}

When $X$ is a smooth Fano surface over a number field, one can prove a
little more:

\begin{theorem}\label{fanotheorem}
Let $X$ be a smooth, projective Fano surface, defined over a number
field $K$.  Then for all but finitely many primes $p$, $X$ has
ordinary reduction at $p$ and moreover the de Rham-Witt cohomology of
$X_p$ is torsion free.
\end{theorem}

\begin{proof}
     It follows from the results of \cite{mehta97} and \cite{hara96},
that if $X$ is a smooth, projective and Fano variety $X$ over a number
field, for all but finitely many primes $p$, the reduction $X_{p}$ is
$F$-split.  By Theorem~\ref{surface-fsord} we see the reduction modulo
all but finitely many primes $p$ of $K$ gives an ordinary surface.
Then by Lemma~9.5 of \cite{bloch86} and Proposition~\ref{fintor}, the
result follows.
\end{proof}

\begin{example} \label{fano} Let $K=\Q$ and
$X\subset \P^n$ be any Fermat hypersurface of degree $m$ and $n\geq
6$. If $m<n+1$ then this hypersurface is Fano but by \cite{toki96}
this hypersurface does not have Hodge-Witt reduction at primes $p$
satisfying $p\not\equiv 1\mod m$. This gives examples of Fano
varieties which are ($F$-split but are) not Hodge-Witt or ordinary.
\end{example}

\begin{remark}
It is clear from Example~\ref{fano} that there exist Fano varieties
over number fields which do not have Hodge-Witt reduction modulo an
infinite set of primes and thus this indicates that in higher
dimension $p_g(X)$ is not a good invariant for measuring this
behavior. The following question and subsequent examples suggests
that the Hodge level may intervene in higher dimensions.
\end{remark}

\begin{question}\label{hodge-level}
   Let $X/K$ be a smooth, projective Fano variety over a number
field. Assume that $X$ has Hodge level $\leq 1$ in the sense of
\cite{deligne72}. Then does $X$ have Hodge-Witt reduction modulo  all
but a finite number of primes of $K$?
\end{question}

\begin{remark}
   A list of all the smooth complete intersection in $\P^n$ which
are of Hodge level $\leq 1$ is given in \cite{rapoport72} and one
knows from \cite{suwa93} that complete intersections of Hodge level 1
are Hodge-Witt.
\end{remark}

\begin{remark}
The first author has answered the Question~\ref{hodge-level}
affirmatively (in \cite{joshi00b}) for $\dim(X)=3$ where the Hodge
level condition is automatic.
\end{remark}

   Recall that an abelian variety $A$ over a perfect field is
Hodge-Witt if and only if the $p$-rank of $A$ is at least $\dim(A)-1$
(see \cite{illusie83a}).  This together with
Theorem~\ref{positive-density} gives

\begin{theorem}
   Let $A/K$ be an abelian threefold over a number field
$K$. Then there exists a set of primes of positive density in $K$ such
that $A$ has Hodge-Witt reduction at these primes.
\end{theorem}

\subsection{Hodge-Witt torsion}
We include here some  observations
probably well-known to the experts, but we have not
found them in print.  We assume as in
the previous section that $X/K$ is smooth projective variety over a
number field and that we have fixed a regular, proper model smooth
over some open subscheme of the ring of integers of $K$ and whose
generic fiber is $X$.

Before we proceed we record the following:

\begin{proposition} Let $X/K$ be a smooth projective variety over a
number field $K$.  Then there exists an integer $N$ such that for all
primes $\wp$ in $K$ lying over any rational prime $p\geq N$, the
following dichotomy holds 
\begin{enumerate} 
\item either for all
$i,j\geq 0$, the Hodge-Witt groups are free $W$-modules (of finite
type), or 
\item there is some pair $i,j$ such that
$H^i(X_\wp,W\Omega_{X_\wp})$ has infinite torsion.  
\end{enumerate}
\end{proposition} 

\begin{proof} 
   Choose a finite set of primes $S$ of $K$, such that $X$ has a
proper, regular model  over $\Spec(\O_{K,S})$, where $\O_{K,S}$
denotes the ring of $S$-integers in $K$. Choose $N$ large enough
so that for any prime $\wp$ lying over a rational prime $p>N $,
we have $\wp \not\in S $, and all  the crystalline cohomology
groups of $X_\wp$ are torsion-free. We note that this choice of
$N$ may depend on the choice of a regular proper model for $X$
over $\Spec(\O_{K,S})$. If $\wp$ is such that
$H^i(X_\wp,W\Omega^j_{X_\wp})$ are all finite type,  then by the
degeneration of the slope spectral sequence at the $E_1$-stage by
Bloch-Nygaard  (see Theorem~3.7 of \cite{illusie79b}), and the
fact that the crystalline cohomology groups are torsion free, it
follows that the Hodge-Witt groups are free as well. If, on the
other hand,  some Hodge-Witt group of $X_\wp$ is not of finite
type over $W$, then we are in the second case.  \end{proof}

\begin{question}\label{hodge-witt-torsion}
   Let $X/K$ be a smooth projective variety over a number
field.  When does there exist an infinite set of primes of
$K$ such that the Hodge-Witt cohomology groups of the reduction
$X_\wp$ at $\wp$ are not Hodge-Witt? 
\end{question}

We would like to explicate the
information encoded in such a set of primes (when it exists).  

\begin{proposition}
Let $A/K$ be an abelian surface over a number field $K$. Then there
exists infinitely many primes $\wp$ such that
$H^2(X_\wp,W(\O_{X_\wp}))$ has infinite torsion if and only if there
exists infinitely many primes $\wp$ of supersingular reduction for
$X$. In particular,  let $E$ be an elliptic curve over $\Q$ and let
$X=E\times_{\Q}E$. Then for an infinite set of primes of $\Q$, the
Hodge-Witt groups $H^i(X_p,W\Omega_{X_p}^j)$ are not torsion free for
$(i,j)\in\left\{(2,0),(2,1)\right\}$.
\end{proposition}

\begin{proof}
The first part  follows from  the results of
\cite[Section~7.1(a)]{illusie79b}.  The second part follows from
combining the first part  with 
Elkies's theorem (see \cite{elkies87}),  that given an elliptic curve
$E$ over $\Q$,  there are infinitely
many primes $p$ of $\Q$ such that $E$ has supersingular reduction. 
\end{proof}

\begin{example} Results of \cite{toki96} on Fermat hypersurfaces
together with foregoing discussion indicate similar examples as above. 
  These are the only examples of this
phenomena we know so far related to the above question. 
\end{example}


\end{document}